# Comment: Classifier Technology and the Illusion of Progress

**Jerome H. Friedman**

This paper provides a valuable service by asking us to reflect on recent developments in classification methodology to ascertain how far we have progressed and what remains to be done. The suggestion in the paper is that the field has advanced very little over the past ten or so years in spite of all of the excitement to the contrary.

It is of course natural to become overenthusiastic about new methods. Academic disciplines are as susceptible to fads as any other endeavor. Statistics and machine learning are not exempt from this phenomenon. Often a new method is heavily championed by its developer(s) as the "magic bullet" that renders past methodology obsolete. Sometimes these arguments are accompanied by nontechnical metaphors such as brain biology, natural selection and human reasoning. The developers become gurus of a movement that eventually attracts disciples who in turn spread the word that a new dawn has emerged. All of this enthusiasm is infectious and the new method is adopted by practitioners who often uncritically assume that they are realizing benefits not afforded by previous methodology. Eventually realism sets in as the limitations of the newer methods emerge and they are placed in proper perspective.

Such realism is often not immediately welcomed. Suggesting that an exciting new method may not bring as great an improvement as initially envisioned or that it may simply be a variation of existing methodology expressed in new vocabulary often elicits a strong reaction. Thus, the messengers who bring this news tend to be, at least initially, unpopular among their colleagues in the field. It therefore takes courage to provide this type of service, and Professor Hand is to be congratulated for this thoughtful article.

Of course, simply because new methodologies are often overhyped does not necessarily imply that they do not, at least sometimes, represent important progress. In the case of classification, I believe that there have been major developments over the past ten years that have substantially advanced the field, both in terms of theory and practice. Although I find myself in agreement with most of the premises of this article, I do not see how they lead to the implication that such advances are "largely illusionary."

There appear to be three main premises presented in the article. First, the improvements realized by the newer methods over the previous ones are less than those achieved by the previous ones over their predecessors, presumably no methodology at all. Second, the evidence often presented (at least initially) in favor of the superiority of the newer methods is often suspect. Finally, the newer methods do not solve all of the outstanding important problems that remain in the field of classification. In my view these observations are correct and underappreciated in the field. The article does an important service by illustrating them so forcefully. However, the truth of these assertions does not imply lack of important progress; only that low-lying fruit is often easier to gather, we should be more thorough concerning validation when initially presenting new procedures and there is still important work to be done.

One of the main assertions in the paper is that, in many applications, older methods often yield error rates comparable to the more modern ones. This is of course true and is intrinsic to the classification problem, especially when the metric used to measure performance is based on error rate. First, there is the irreducible error caused by the fact that the predictor variables $\mathbf{x}$ often do not contain enough information to specify a unique value for the outcome variable $y$. At best, they specify a probability distribution of possible values $\Pr(y|\mathbf{x})$ which is hopefully different for differing values of $\mathbf{x}$, indicating some predictive power. This phenomenon afflicts all prediction problems. A second phenomenon is peculiar to classification; it is not necessary to accurately estimate $\Pr(y|\mathbf{x})$ to achieve minimal error rate. All that is required of the estimates $\widehat{\Pr}(y|\mathbf{x})$ is

$$\tag{1} \arg\max_y \widehat{\Pr}(y|\mathbf{x}) = \arg\max_y \Pr(y|\mathbf{x}).$$

*Department of Statistics and Stanford Linear Accelerator Center, Stanford University, Stanford, California 94305, USA (e-mail: jhf@stanford.edu)*

The actual values of the estimates for differing values of $y$ need not be close to their respective underlying true values. The estimates for the nonmaximizing probabilities need not even be in the correct order. Thus, more flexible (modern) procedures that are better able to estimate more complex probability structures need not produce dramatically lower error rates in many applications. This also accounts for the "flat minimum" effect discussed in the paper.

As pointed out in the paper, classification procedures are often used in contexts where error rate is not the relevant quantity; functionals of $\Pr(y|\mathbf{x})$ other than (1) are of interest. For example, in many two-class classification problems $y \in \{-1, 1\}$, the important quantity is the rank order of $\{\Pr(y = 1|\mathbf{x}_i)\}_{i \in T}$, where $T$ is a set of observations with unknown outcome. In other applications, interest is in the actual probabilities themselves. In such settings it is likely that more accurate estimates of $\Pr(y|\mathbf{x})$ afforded by more flexible modern techniques will yield distinctly superior results to the older less flexible methods, even though their respective error rates are not dramatically different. The paper properly criticizes the classification literature for presenting comparisons mostly in terms of error rate, even though this is the criterion used for nearly all of the classification comparisons presented in the paper.

The primary evidence intended to suggest lack of progress is the comparisons presented in Table 1. Here the error rate of an older method, linear discriminant analysis (LDA), is compared with that of the current best method for each of a selected set of problems. In spite of the general insensitivity of error rate as differentiating criterion (as discussed above), LDA seems to produce distinctly inferior results in many of these problems. In more than half of the examples, its error rate is at least 45% greater; in one example, it is nearly six times as great. Of course there is a selection bias of unknown magnitude in choosing the best method, but it is difficult to conclude from the evidence presented that LDA is competitive with the best current methods, even in terms of error rate. The paper suggests that large ratios in small error rates "will correspond to only a small proportion of new data points." This is true but not relevant. If a zip code classifier makes twice as many errors, it costs the post office twice as much to handle the misdirected mail. I have yet to see a problem where costs are proportional to the Prop linear statistic shown in the last column of Table 1.

The paper presents a regression example (Section 2.1) to illustrate that including additional predictor variables that are highly correlated with those that are already part of the analysis produces little gain in performance. This is true of all methods, old and new, and no evidence is provided to suggest that older methods are better able to incorporate additional information from such variables.

A second principal premise of the paper is that the evidence for the superiority of new methods is generally based on empirical comparisons which are susceptible to major weaknesses that place their validity in question. I could not be in more agreement with this point. Section 5 of the paper should be required reading for all practitioners and researchers in the field. In my data mining course, I have a lecture called "comparison caution" that addresses many of the same issues. Empirical comparisons should be viewed with skepticism, especially when the authors' new method is one of the competitors. Even when this is not the case, the authors performing the study often have a favorite technique which usually emerges as the top performer. When interpreting such studies, I tend to ignore the apparent top performer and look at the relative rankings of the other methods, presuming that the authors have less expertise and vested interest in them. Even when a comparison is free of all of the biases discussed in Section 5, its results should not be extrapolated beyond the specifics of the problem represented by the data set being used. All methods have particular problems for which they are especially well suited and others for which they are not. Sometimes only a minor change in the problem setup can produce substantial changes in performance rankings. Results of empirical comparisons can be useful, especially when aggregated over time, but the natural tendency to overinterpret individual studies should be avoided. Of course, the same caution should be applied to the empirical comparisons presented in this paper.

Simply because the initial evidence for the superiority of a method can be questioned does not necessarily imply that is not useful or that it does not represent progress. Practitioners try various methods and, as time evolves, some emerge as being more useful that others. Many of the "new" proposals of the distant past have not survived the test of time and are now long forgotten. Those that have emerged as being generally useful, such as logistic regression, LDA and decision trees, have survived to see common use. No one is claiming that all of the new techniques proposed in the literature over past ten years represent major advances. However, I believe that a body of evidence is emerging that suggests that some of them, such as the ensemble methods (bagging and boosting) and support vector machines, offer substantial advantages over the earlier methods in enough situations to be regarded as major advances. This is especially the case in scientific and engineering applications, where decision boundaries are often complex and far from being linear.

Another major premise of the paper is that there are important issues that affect classification performance that are not addressed by most modern methodology. These include population drift, sample selectivity bias, errors in class labels and arbitrariness in class definitions. Again I could not agree more. Issues of nonrepresentative training data tend to be overlooked by the academic community, although they are probably well known to most practitioners. (See [3]. I spend several lectures in my data mining course covering these topics.) Obtaining high-quality representative training data is generally more important to success than choice of a particular classifier, although given such data, choosing the best classifier can often provide considerable additional benefit. In many data mining applications, the data were collected for a different purpose than solving the current problem and one does not have influence over its quality or value. The analyst is forced to do the best that can be done with the data at hand.

The problem of training data being different from future data to be predicted is common to all prediction, not just classification. The fundamental issue is similar whether the differences arise through random sampling from a static population or are caused by one of the more deterministic mechanisms cited in the paper. As noted in the paper, the antidote is to limit reliance on the training data by not fitting it too closely. This is the basic principal underlying regularization. The paper argues that older methods are "simpler," thereby inducing more regularization, which in turn causes them to be more resistant to these types of problems. This need not be the case.

Almost all of the modern procedures incorporate a regularization parameter that controls the degree to which they are allowed to fit the training data. By adjusting the value of this parameter, one can produce a sequence of models of increasing complexity from the very simplest that makes the same prediction everywhere to highly complex functions that capture the fine details of the predictive relationship as reflected in the training data. Highly regularized versions of different procedures may capture somewhat different aspects of the gross features of the probability distribution, but in the absence of knowledge concerning the nature of the population drift, there is no a priori reason to suspect that one is better than the other. An important consequence of the presence of population drift and related problems is that model selection based on traditional techniques such as bootstrapping or cross-validation becomes overly optimistic; they will tend to produce insufficient regularization. Thus, care must be taken to regularize more heavily than suggested by these model selection techniques when such problems are suspected.

Most older classification methods limit the degree to which one can control the amount of regularization. It it not clear that the amount arbitrarily applied by these procedures is necessarily appropriate in any particular problem. In fact there are many situations in which older methods provide insufficient regularization. This is especially the case in modern analytical chemistry and bioinformatics applications, where there are many more predictor variables than training observations, and simple logistic regression and LDA completely fail. There has been considerable recent research that has led to modern classification methods than allow the application of more regularization than the older traditional methods. These, in my view, also represent major progress.

Errors in class labels is a classic robustness issue. Estimation in the presence of badly measured outcomes has been extensively studied in the regression literature, but less so in classification. As in regression, the solution is to employ loss criteria that are less sensitive to individual extreme measurements. It has been suggested that logistic likelihood and the support vector machine hinge loss are more robust to misspecification of class labels than squared-error loss or, especially, the exponential loss associated with AdaBoost, since they weight realized outcomes of low estimated probability less heavily. Even more robust (nonconvex) loss criteria have been proposed for classification (see [1, 2]). Some older methods such as logistic regression should be fairly robust to mislabeling, but others like LDA are likely to exhibit poor robustness properties; estimates of the pooled covariance matrix can be highly distorted by only a few mislabeled observations, especially at the extremes of the data distribution.

The problem of arbitrariness of class labels is often caused by trying to make the problem conform to the method rather than the other way around. If an outcome variable realizes continuous numeric values, then it should be treated as such and regression rather than classification technology would be more appropriate. There have been recent important advances in regression methodology that parallel those in classification. If thresholding numeric variables to create a classification problem happens to be appropriate and the class labels have changed, then, as the paper suggests, one can simply retrain the classifier with the new definitions. This requires that the original raw data be saved. Given the very low cost of storage media, this should always be encouraged for a wide variety of reasons.

Recent research has not solved all of the outstanding problems in the field of classification, especially those associated with nonrepresentative training data. All procedures are vulnerable to these effects and, as discussed above, it is not clear that the older methods enjoy more immunity than the more recent ones. Also, these problems are more prevalent in the commercial sector involving financial and consumer behavior applications than in scientific and engineering fields where the laws governing the systems under study tend to be more stable. Nevertheless, solutions to these problems would also represent major advances. The paper does an important service by directing our attention to them, but this does not imply that there has not been substantial progress in other important aspects of the classification problem in the recent past.

Whether or not a new method represents important progress is, at least initially, a value judgement upon which people can agree to disagree. Initial hype can be misleading and only with the passage of time can such controversies be resolved. It may well be too soon to draw conclusions concerning the precise value of recent developments, but to conclude that they represent very little progress is at best premature and, in my view, contrary to present evidence.

I thank Professor Hand for this thoughtfully provocative article. It gives all of us an opportunity to look past our enthusiasm and take a deeper look at the remaining central issues. I look forward to research that produces solutions to these outstanding problems and to future discussions as to whether they represent major progress. Finally, I would like to add another relevant quote to that of Eric Hoffer mentioned in the article. This one is attributed to Yogi Berra: "Prediction is difficult, especially when it's about the future."